\documentclass[a4paper,11pt]{article}
\usepackage[dvips]{graphicx}
\usepackage{amssymb}
\def\intm{{\,\int\!\!\!\!\!\!\!-\!\!\!-}}
\def\intms{{\,\int\!\!\!\!\!-\,}}

\newtheorem{theo}{Theorem}[section]
\newtheorem{lem}[theo]{Lemma}
\newtheorem{pro}[theo]{Proposition}
\newtheorem{rem}[theo]{Remark}
\def\R{\mathbb{R}}
\title{A Gamma-convergence argument for the blow-up\\
of a non-local semilinear parabolic equation with Neumann boundary conditions}
\author{A. El Soufi, M. Jazar\footnote{The second author is supported by a grant from
Lebanese University.}\, and R. Monneau\\ Universit\'e de Tours,
Lebanese University and ENPC}

\begin{document}

\maketitle

\noindent{\small{\textbf{ Abstract.} In this paper we study a
simple non-local semilinear parabolic equation with Neumann
boundary condition. We give local existence result and prove
global existence for small initial data. A natural non increasing
in time energy is associated to this equation. We prove that the
solution blows up at finite time $T$ if and only if its energy is
negative at some time before $T$. The proof of this result is
based on a Gamma-convergence technique.}}

\noindent{\small\textbf{AMS Subject Classifications:}}
{\small{Primary 35B35, 35B40, 35K55; secondary 35K57,
35K60}}\hfill

\noindent{\small\textbf{Keywords:}} {\small{ Semi-linear parabolic
equations, Blow-up, Global existence, Asymptotic behavior of
solutions, Gamma convergence, Neumann Heat kernel, Non local term,
Comparison principle}}

\section{Introduction}

\subsection{Setting of the problem}

In this paper, we consider a bounded domain $\Omega$ in $\R^N$
that is uniformly regular of class $C^2$, and we study the
solutions $u(x,t)$ of the following equation for some $p>1$
(denoting the mean value $\frac{1}{|\Omega|}\int_\Omega f$ by
$\intms_\Omega f$ for a general function $f$):
\begin{equation}\label{eq:1}
\left\{\begin{array}{lcl}
\displaystyle{u_t-\Delta u = |u|^p-  \intm_\Omega |u|^p \quad \mbox{on}\quad \Omega} \\
\\
\displaystyle{\frac{\partial u}{\partial n}=0 \quad \mbox{on}\quad
\partial\Omega}
\end{array}\right.
\end{equation}
with the initial condition
\begin{equation}\label{eq:2}
\left\{\begin{array}{l}
\displaystyle{u(x,0)=u_0(x) \quad \mbox{on}\quad \Omega}   \\
\displaystyle{\mbox{with} \quad \intm_\Omega u_0=0}
\end{array}
\right.
\end{equation}
It is immediate to check that the integral (or the mean value) of
$u$ is conserved (at least once you precise the meaning of the solution). \\
Stationary solutions of Equation (\ref{eq:1}) are in fact
critical points of the energy functional
$$E(u)=\int_\Omega \left[ \frac12 |\nabla u|^2 -\frac1{p+1}u|u|^p\right]$$
under the constraint that $\intms_\Omega u$ is equal to a given constant.\\
Without loss of generality, we can assume that $|\Omega|=1$.
Indeed, if $u$ is a solution of (\ref{eq:1})-(\ref{eq:2}) in
$\Omega$ and $\lambda>0$, then $v(t,x):=\lambda^{\frac2{p-1}}
(\lambda^2 t,\lambda x)$ is a solution in $\lambda^{-1}\Omega$.
Throughout the paper, we assume $|\Omega|=1$, except when the
volume $|\Omega|$ is explicitly mentioned to show the dependence
of the constants.

\subsection{Motivation of the problem}

A lot of work has been done on scalar semilinear parabolic
equations whose the most famous example is
$$u_t-\Delta u =u^p$$
and the problem of global existence or blow-up is quite well
understood (see for instance \cite{Fujita,Ball} for an energy
criterion for blow-up, \cite{GK} for a study of self-similar
blow-ups; see also \cite{Zaag,S3} and the numerous references
therein). Of course, the Maximum Principle plays a fundamental
role in the establishment of results in this setting.
\textbf{However, concerning the problem of describing the blow-up
set, very few results are known. For instance, in dimension 2, the
question of whether there exists a solution whose blow-up set is
an ellipse is still unanswered. Recently, Zaag \cite{Zaag1}
established the first regularity results for the blow-up set,
based on global estimates independent of the blow-up point
obtained by Merle and Zaag \cite{MZ} through the proof of a
Liouville theorem.}

In the case of parabolic systems or non-local scalar parabolic
equations, even if some Maximum Principles may hold, it is often
necessary to introduce new techniques. One of the most famous
examples is the Navier-Stokes equation (see \cite{L}), which can
be written on the vorticity $\omega =\mbox{ curl } u$ (with $u$
the velocity and $e=\frac12(\nabla u+{}^t \nabla u)$ the
deformation velocity):
$$\omega_t-\nu \Delta \omega=-(u \cdot \nabla ) \omega +e\cdot \omega$$
where the right hand side is non-local and quadratic in $\omega$.
If we consider this equation on $\Omega =\R^3\backslash {\mathbb
Z}^3$, the
following quantity is conserved by the equation
$$\int_\Omega \omega(t,x)\ dx .$$

One of the simplest examples of non-local and quadratic equation
is
$$u_t-\Delta u = u^2-\frac{1}{|\Omega|}\int_\Omega u^2$$
with Neumann boundary condition on $\partial\Omega$ so that the
quantity $\int_\Omega u(t,x)\ dx$ is conserved. This equation is
also related to Navier-Stokes equations on an infinite slab for
other reasons explained in \cite{BDS2}. Problem
(\ref{eq:1})-(\ref{eq:2}) is a natural generalization of this
latter for which we provide a global existence result for small
initial data as well as a new blow-up criterion based on partial
Maximum Principles and on a Gamma-convergence argument.

\subsection{Main results}

In this subsection, we present our main results: local existence,
global existence for small initial data, energetic criterion for
blow-up of solutions based on an optimization result  of
independent interest that we prove by a Gamma-convergence
technique. Furthermore we give a global existence result in the
case $p=2$, expliciting the constants as a function of the
geometry
of the domain.\\

First, let us mention that the classical semigroup theory enables
us to prove, more or less directly:
\begin{itemize}
\item local existence and uniqueness of solutions of
(\ref{eq:1})-(\ref{eq:2}) for any initial data $u_0\in
C(\overline{\Omega})$ (see the next section),

\item global existence and exponential decay of solutions of
(\ref{eq:1})-(\ref{eq:2}) for small initial data $u_0\in
C(\overline{\Omega})$. That is, there exists some (implicit)
constant $\rho>0$ depending on the geometry of $\Omega$, so that
$\|u_0\|_{L^\infty}<\rho$ implies global existence and exponential
decay of $u$: $\|u(t)\|_{L^\infty}\le Ce^{-\alpha t}$ for some
positive constants $C$ and $\alpha$.
\end{itemize}

To prove results of this kind it suffices for instance to follow
the arguments of the proof of \cite[Proposition 5.3.9]{CH}. See
also \cite{BDS} for a 1-dimensional result in this direction.\\

Our main purpose is to give a natural energetic criterion for the
blow-up in finite time of solutions of (\ref{eq:1})-(\ref{eq:2})
in the case $1<p\leq 2$. \textbf{Our proof relies on the same main
idea introduced by Levine \cite{Levine} and Ball \cite{Ball} in
the sense that the blow-up will follow from a nonlinear
differential inequality that we show to be satisfied by the
$L^2$-norm of the solution.}

First, it is quite easy to see that, $\forall p>1$, the energy
\begin{equation}\label{eq:E}
E(u)=\int_\Omega \frac12 |\nabla u|^2 -\frac1{p+1}u|u|^p
\end{equation}
of a solution $u$ of (\ref{eq:1})-(\ref{eq:2}) is non increasing
in time (see Proposition \ref{pro:1}). Our main result in this
direction is the following

\begin{theo}\label{th:4}\textbf{ (Energetic criterion for blow-up, case
$1<p\leq 2$)}\\
Let us assume that $p\in (1,2]$ and let $u$ be a solution of
(\ref{eq:1})-(\ref{eq:2}) with $u_0\in C(\overline{\Omega})$,
$u_0\not\equiv 0$. If the energy of $u_0$ is nonpositive, that is
$E(u_0)\leq 0$, \textbf{then the solution does not exist in
$L^2(\Omega)$ for all $t>0$. Moreover, there exists $T>0$ such
that if $u\in L^\infty_{loc}([0,T);L^2(\Omega))$, then}
\begin{equation}\label{eq:0}
\lim_{t\to T^-}||u(t)||_{L^2(\Omega)}=+\infty .
\end{equation}
\end{theo}

\begin{rem}\label{rem12}
Note that Theorem \ref{th:4} does not imply that the $L^2$-norm of
$u(t)$ blows-up in finite time. Indeed, the solution may simply
not exist till time $T$.

Recall that, for the semilinear heat equation $u_t = \Delta u +
u^p$ on a bounded domain, the generic blow-up profile is given by
(see \cite{HV} for details)
$$u(x,t)\rightarrow u^*(x)\qquad\mbox{as }t\rightarrow T,$$
with
$$u^*(x)\sim C(p){\left[\frac{|\log|x||}{|x|^2}  \right]}^{\frac1{p-1}}$$
$C(p)$ being a constant. Hence, the $L^2$-norm stays generically
bounded whenever $p
> 1+\frac4N$, and blows-up when $p < 1+\frac4N$.
\end{rem}

The condition of nonpositivity of the energy in Theorem \ref{th:4}
is also necessary in the sense that, if the $L^2$-norm of $u(t)$
blows-up at a time $T>0$, then the energy $E(u(t))$ needs to be
negative at some time $t<T$. The situation is even worse: the
energy $E(u(t))$ needs to blow-up to $-\infty$ at a time $T'\leq
T$. Moreover, this property is valid for any $p\in (1,+\infty)$.
Indeed, we have the following

\begin{theo}\label{th:5}\textbf{ ($L^2$ bound on $u$ for bounded from below
energy, $p\in (1,+\infty)$)}\\
Let $p$ be any real number in $(1,+\infty)$ and let $u$ be a
solution of (\ref{eq:1})-(\ref{eq:2}) with $u_0\in
C(\overline{\Omega})$. If there exists a constant $C_0>0$ and a
time $T_0>0$ such that
$$E(t)\geq -C_0 \quad \mbox{for}\quad t\in [0,T_0),$$
then there exists a constant $C>0$ such that
$$||u(t)||_{L^2(\Omega)}\leq C  \quad \mbox{for}\quad t\in [0,T_0).$$
\end{theo}

The case $p>2$ is still not completely understood for us, however,
we believe that the blow-up phenomenon of Theorem \ref{th:4}
occurs for any $p$ in $(1,+\infty)$ and formulate the following
conjecture:

\bigskip

\noindent\textbf{Conjecture} \textbf{ (Energetic criterion for
blow-up, case $\mathbf{p>2}$)}\\
{\it For $p>2$, we conjecture that if $u$ is a solution of
(\ref{eq:1})-(\ref{eq:2}) with $E(u_0)\leq 0$ and $u_0\not\equiv
0$, then $u(t)$ blows-up in finite time.}

\bigskip

The proof of Theorem \ref{th:4} is based, on one hand, on the use
of maximum principles, and, on the other hand, on the following
estimate of independent interest, proved by gamma-convergence:
\begin{theo}\label{th:optimization}\textbf{ (Optimization under a $L^2$ constraint)}\\
For $p>1$, there exists $\theta_0>0$ and $C>0$ such that
\begin{equation}\label{inf}
\inf_{v\in A}\int_\Omega v|v|^{p}+\theta |\nabla v|^2\ge
C\sqrt{\theta}
\end{equation}
for all $\theta\in [0,\theta_0]$, where
$$\displaystyle{A:=\left\{v\in L^{p+1}(\Omega)\; \Big| \; \int_\Omega
v=0, \; \int_\Omega v^2=1 \; \mbox{and} \; v\ge -1\ \mbox{on}\
\Omega\right\}}.$$
\end{theo}

Let us mention that the profile of blowing-up solutions for this
equation seems to us an open problem in general. Besides an
example given in \cite{BDS} of a profile of a blowing-up solution
whose the positive part concentrates at one point in the
one-dimensional case, we do not know if it is possible to build
blowing-up solutions with different profiles.\\


Our next purpose is to focus on global existence results in the
case $p=2$. As mentioned previously, in usual global existence
results, the constant $\rho$ that determines the smallness of the
initial data should depend on the geometry of the domain $\Omega$.
It is interesting to understand this dependence. That is precisely
the aim of our Theorem \ref{th:global}. In particular, we relax
here the assumption $|\Omega|=1$ to show the dependence on the
volume $|\Omega|$. We need first to introduce the following two
invariants:
\begin{itemize}
\item the first positive eigenvalue $\lambda_1(\Omega)$ of the
Laplacian in $\Omega$ under Neumann boundary conditions. Recall
that we have the following isoperimetric-type inequality due to
Szeg\"o \cite{Sz} and Weinberger \cite{Wein}:
\begin{equation}\label{iso}\lambda_1(\Omega)|\Omega|^{2/N}
\leq \lambda_1^*(N),\end{equation} where $\lambda_1^*(N)$ is the
first positive Neumann eigenvalue of the $N$-dim-ensional
Euclidean ball of volume 1.

\item the constant $H(\Omega)$ defined as the supremum over
$(0,+\infty)\times \Omega$ of the function
$t^{N/2}[K(t,x,x)-{1\over |\Omega| }]$ where $K(t,x,y)$ is the
heat kernel associated to the Laplacian in $\Omega$ with Neumann
boundary conditions (see \cite{Davies} and section 2.2 for the
precise definition of $K$ and the existence of $H(\Omega)$).
Notice that one has (see remark \ref{rem:H})
$$H(\Omega)\geq (4\pi)^{-\frac N2}.$$
\end{itemize}

It is also well known that the constant $H(\Omega)$ is closely
related to the so-called Neumann Sobolev constant $C(\Omega)$
defined as the best constant in the inequality : $\forall f\in
H^1(\Omega)$ such that $\int_\Omega f =0$, we have
$\|f\|^2_{2N\over{N-2}} \le C(\Omega) \|\nabla f\|^2_2$ (see, for
instance, \cite[section
3]{JW} for results about this relationship).\\

Let us first remark that we have the following property for $p=2$
$$\mbox{if}\quad ||u_0||_{L^\infty (\Omega)}\leq
\frac{3}{2}\lambda_1(\Omega),\quad \mbox{then}\quad E(u_0)\geq
0,$$ which may indicate (from Theorem \ref{th:4}) that the
corresponding solutions may not necessary  blow-up in finite time.
This shows in particular that it is natural to compare
$||u_0||_{L^\infty (\Omega)}$ with the first eigenvalue
$\lambda_1(\Omega)$ as it can also be seen from the scaling of the
equation for $p=2$.\\

The following theorem gives a global existence result under an
explicit smallness condition on the initial data, depending on
$\lambda_1(\Omega)$ on the one hand and on $N$, and $H(\Omega)$ on
the other hand. For simplicity, we state it only for $p=2$
although a general version is possible.

\begin{theo}\label{th:global}\textbf{ (Global existence for small initial
data with explicit constants, case $p=2$)}. Let $\rho(\Omega)$
be the constant given by
$$\displaystyle{\frac{\rho(\Omega)}{\lambda_1(\Omega)}=\frac{1}{2N}\cdot
\exp\left(-\gamma_{N}\, \lambda_1^*(N)\,
      (H(\Omega))^{\frac{2}{N}}\right)},$$
where $$\gamma_N=\frac{2^7}{N}$$ For every $u_0\in
C(\overline{\Omega})$ satisfying $\int_\Omega u_0=0$ and
\begin{equation}
\|u_0\|_{L^\infty(\Omega)}\le  \rho(\Omega),
\end{equation}
the (unique) solution of the problem (\ref{eq:1})-(\ref{eq:2}) is
defined for all $t$ and tends to zero as $t\to\infty$.
\end{theo}

\begin{rem}
We also provide an exponential decay (see Theorem \ref{th:globalbis}).
\end{rem}

Note that $H(\Omega)$ is invariant by dilation of the domain $\Omega$, and then that
$\rho(\lambda \Omega) = \lambda^{-2}\rho(\Omega) <\rho(\Omega)$ for $\lambda
>1$, but there is no reason in general to get $\rho(\Omega)\le
\rho(\Omega')$ if $\Omega' \subset \Omega$.

\begin{rem}\label{}
Actually, the volume of $\Omega$ being fixed, one could reasonably
expect that the constant $\rho(\Omega)$ is maximal when $\Omega$
is a ball.
\end{rem}


\subsection{Brief review of the literature}

Parabolic problems involving non local terms have been recently
studied extensively in the literature (see for instance
\cite{E,FG,RS}). For local existence and continuation results for
general semilinear equations under the Neumann boundary condition
setting, one can see for example \cite{A} and \cite{Stewart}. In
\cite[Appendix A]{S1} Souplet gives very general local existence
results for a large class of non local problems in time and in
space but in the Dirichlet boundary condition setting. The problem
treated in the present work has been first considered by Budd,
Dold and Stuart (\cite{BDS}) for $p=2$ and in the one dimensional
case. They obtained a theorem like our Theorem \ref{th:global} as
well as a blow-up type result for solutions whose Fourier
coefficients of the initial data satisfy an infinite number of conditions.

Hu and Yin (\cite{HY}) considered slightly different problems. They
showed in particular blow-up result (see \cite[Theorem 2.1]{HY}),
based on energy criteria, considering $u|u|^{p-1}$ instead of
$|u|^p$. They showed also (see \cite[Theorem 3.1 and 3.2]{HY})
global existence for positive solutions and $p$ not too large. A
radial blowing-up solution for $p$ large is also given.

Wang and Wang (\cite{Wang}) considered a more general problem of
the form
\begin{equation}\label{wang}u_t-\Delta u=ku^p-\int
u^q\end{equation} with Neumann or Dirichlet boundary conditions
and positive initial data. They showed global existence and
exponential decay in the case where $p=q$, $|\Omega|\le k$ and
Neumann boundary condition. They also obtain a blow-up result
under the assumption that the initial data is bigger than some
"gaussian function" in the case where $|\Omega|>k$.

Finally, in \cite[Theorem 2.2]{S2}, Souplet determines exact
behavior of the blow-up rate for equations of the form
(\ref{wang}) with $k=1$ and $p\not= q$.

\subsection{Organization of the article}
Our paper is organized as follows. In the second section we first
set the space under which problem (\ref{eq:1})-(\ref{eq:2}) admits
a unique local solution. Section 3 is devoted to the proof of
Theorems \ref{th:4} and \ref{th:5}. In section 4 we give the proof
of the optimization result (Theorem \ref{th:optimization}) which
is based on a result of Gamma-convergence of Modica \cite{modica}.
For the convenience of the reader we provide in the appendix
(section 6)  a self-contained proof of the corresponding
Gamma-convergence-like result. In section 5, we give the proof of
Theorem \ref{th:global}.

\section{Local existence result}

We recall that $\Omega$ is bounded. The basic space to be
considered in this paper is the space $C(\overline{\Omega})$ of
continuous functions. Following the notations \footnote{Since
$\Omega$ is bounded and $\partial\Omega\in C^2$, we can
  easily check that $D^q$ is dense in $C(\overline{\Omega})$}
of Stewart \cite{Stewart} denote for $q>N$ by
$$D^q:=\left\{u\in C(\overline{\Omega});\,\,u\in W^{2,q}(\Omega),\,
\Delta u\in C(\overline{\Omega}),\,\mbox{and}\,\frac{\partial
u}{\partial n}=0 \,\mbox{on}\,\partial\Omega\right\}.$$ Set
$$D:=\bigcup_{N<q<+\infty} D^q.$$
Then we have as a direct application of \cite[Theorem 2]{Stewart}:

\begin{theo} The operator $-\Delta$ with domain $D$ generates
an analytic semigroup in the space $C(\overline{\Omega})$ with the
supremum norm.
\end{theo}

See Lunardi \cite{Lun} for the definition of analytic semigroups.

Then we have
\begin{theo}\label{th:2.2}\textbf{ (Local existence result, $1<p<+\infty$)}\\
For every $u_0\in C(\overline{\Omega})$ there is a $0<t_{max}\le
\infty$ such that the problem (\ref{eq:1})-(\ref{eq:2}) has a
unique mild solution, i.e. a unique solution
$$u\in C\left([0,t_{max});C(\overline{\Omega})\right)\cap
C^1\left((0,t_{max});C(\overline{\Omega})\right)\cap
C\left((0,t_{max});D\right)$$ of the following integral equation
$$u(t)=e^{t\Delta}u_0+\int_0^te^{(t-s)\Delta}f(u(s))\,ds$$
on $[0,t_{max})$, with
$$f(u(s))=|u(s)|^p-\intm_\Omega |u(s)|^p.$$
Moreover, we have
\begin{equation}\label{eq:moyenne}
\int_\Omega u(t) = 0 \quad \mbox{for all}\quad t\in [0,t_{max})
\end{equation}
and if $t_{max}<\infty$ then
$$\lim_{t\uparrow t_{max}}\|u(t)\|_{L^\infty(\Omega)}=\infty.$$
\end{theo}

\noindent\textbf{Proof of Theorem \ref{th:2.2}}\\
First let us remark that the non linearity in (\ref{eq:1}): $u\in
C(\overline{\Omega}) \longmapsto f(u)=|u|^p-\intms_\Omega |u|^p
\in C(\overline{\Omega})$ satisfies the hypothesis of
\cite[Theorem 6.1.4]{Pazy}, namely, $f$ is locally Lipschitz in
$u$, uniformly in $t$ on bounded intervals of time.\\
Then a standard semigroups result (\cite[Theorem
6.1.4]{Pazy}) gives the local existence result. Reminder to
show (\ref{eq:moyenne}). This comes from the simple computation
for $t\in (0,t_{max})$:
$$\frac{d}{dt}\left(\int_\Omega u \right)=\int_\Omega u_t = \int_\Omega
\Delta u + \int_\Omega f(u)= 0$$ because of the definition of $f$,
the integration by part on $\Delta u$, and the Neumann boundary
condition $\frac{\partial u}{\partial n}=0$.
${}\hfill\square$\\

From this result, we see that the solution is a classical solution
of equation (\ref{eq:1}) on $(0,t_{max})\times \Omega$ with
initial condition (\ref{eq:2}), and then from the standard
parabolic estimates (see Lieberman \cite{Lieberman}) and classical
bootstrap arguments, we get
$$u\in C^\infty((0,t_{max})\times \Omega).$$


\section{Blow-up: proofs of Theorems \ref{th:4} and \ref{th:5}}
\textbf{As mentioned in the introduction, we follow the energetic
method introduced by Levine \cite{Levine} and Ball \cite{Ball}.
The main idea is to show that the $L^2$-norm of the solution
satisfies some super-linear differential inequality which implies
the finite time blow-up.}

All along this section, we denote by $u$ a solution of
(\ref{eq:1})-(\ref{eq:2}) whose initial data $u_0$ satisfies
$\int_\Omega u_0 =0$. Also, we will assume, without lack of
generality, that $|\Omega|=1$ so that we have, in particular,
$\intms_\Omega u^2 =\int_\Omega u^2$.

Let us recall the expression of the energy of the problem
$$E(u):=\int\frac 12|\nabla u|^2-\frac 1{p+1}u|u|^p.$$
\begin{pro}\label{pro:1}\textbf{ (Energy decay)}\\
The energy $E(t):=E(u(t))$ is a non increasing function of $t$ in $(0,\infty)$.
\end{pro}
\textbf{Proof of proposition \ref{pro:1}}\\
A direct computation using (\ref{eq:1})-(\ref{eq:2}) and the fact that $\int_\Omega
u_t=0$ yields
$$\frac{d}{dt} E(u(t))=\int_\Omega \left(-\Delta u -|u|^p\right)u_t=-\int_\Omega
(u_t)^2\le 0.$$
$ {}\hfill\square$
\begin{lem}\label{pro:2}\textbf{ ($L^2$ bound from below)}\\
Let us define
$$F(t)=\int_\Omega u^2(t)$$
Then  we have, $\forall p>1$,
$$\frac12 F'(t)= -(p+1)E(t)+ \frac{p-1}{2}\int_\Omega |\nabla u|^2.$$
In particular, we get
$$\frac12 F'(t)\geq  -(p+1)E(0)+ \frac{p-1}{2} \lambda_1(\Omega) F(t).$$
\end{lem}
Consequently, if $E(0)\leq 0$, then $$F(t)\geq F(0)e^{(p-1)\lambda_1(\Omega) t}.$$
\noindent\textbf{Proof of Lemma \ref{pro:2}}\\
The lemma is a consequence of the following computation
$$\begin{array}{ll}
\frac12 F'(t) & =  \displaystyle{\int_\Omega u u_t} \\
\\
& = \displaystyle{\int_\Omega u \left(\Delta u + |u|^p\right)}\\
\\
&=  \displaystyle{\int_\Omega -|\nabla u|^2  + u|u|^p}\\
\\
& = \displaystyle{-(p+1)E(t)+ \frac{p-1}{2}\int_\Omega |\nabla u|^2}
\end{array}$$
${}\hfill\square$

\begin{lem}\label{pro:3}\textbf{($L^2$ bound from below when
$\mathbf{\inf_x u(t)\geq -||u||_{L^2(\Omega)}}$)}\\
Let $0<t_1<t_2<\infty$, $p>1$ and assume that
\begin{equation}\label{eq:10}
\inf_x u(t)\geq -||u||_{L^2(\Omega)} \quad \mbox{for all}\quad
t\in (t_1,t_2).
\end{equation}
Then for all $\beta\in (2,p+1)$, there exists two constants
$C_\beta>0$ and $\lambda_\beta >0$ such that, for all $t\in
(t_1,t_2)$, we have
\begin{equation}\label{eq:11}
||u(t)||_{L^2(\Omega)}\geq \lambda_\beta,
\end{equation}
and
$$ \frac12 F'(t)\geq -\beta E(t)+ C_\beta F(t)^{\frac{p+3}{4}}.$$
\end{lem}

\noindent\textbf{Proof of Lemma \ref{pro:3}}\\
Let us consider a parameter $\beta \in (2,p+1)$. We have
$$\begin{array}{ll}
\displaystyle{\frac12 F'(t)} & =   \displaystyle{\int_\Omega -|\nabla u|^2  + u|u|^p}\\
\\
& = \displaystyle{-\beta \left(\int_\Omega \frac12 |\nabla u|^2 -\frac1{p+1}
  u|u|^p\right) +\frac{\beta -2}{2}\int_\Omega |\nabla u|^2
+\left(1-\frac{\beta}{p+1}\right)\int_\Omega u|u|^p}\\
\\
&=\displaystyle{-\beta E(t)+\frac{p+1-\beta}{p+1}\left(\int_\Omega u|u|^p +\gamma |\nabla
u|^2\right)}
\end{array}$$
with
$$\gamma=\frac{\beta-2}2\times \frac{p+1}{p+1-\beta}$$
Here we will use Theorem \ref{th:optimization}. To this end, we
define
$$\lambda=||u||_{L^2},\quad v=\frac{u}{\lambda}.$$
Then Theorem  \ref{th:optimization} claims that
$$\int_\Omega v|v|^p +\theta |\nabla v|^2\geq C\sqrt{\theta}, \quad
\mbox{if}\quad 0\leq \theta\leq\theta_0.$$ Then we get
$$\begin{array}{ll}
\displaystyle{\frac12 F'(t)} & =  \displaystyle{-\beta
  E(t)+\frac{p+1-\beta}{p+1}\lambda^{p+1}
\left[\int_\Omega v|v|^p +\gamma\lambda^{1-p}|\nabla v|^2\right]}\\
\\
& \geq   \displaystyle{-\beta E(t)+\frac{p+1-\beta}{p+1}
C\gamma^{\frac12}\lambda^{\frac{p+3}{2}}}
\end{array}$$
if $\gamma\lambda^{1-p}\leq \theta_0$. ${}\hfill\square$

\begin{lem}\label{pro:4}\textbf{($\mathbf{L^\infty}$ bound from
below for $\mathbf{1<p\leq 2}$)}\\
Let $p\in (1,2]$, and let $u$ be a solution of the problem
(\ref{eq:1})-(\ref{eq:2}) with the initial data $u(t=0)=u_0$
satisfying the condition
$$u_0\ge -||u_0||_{L^2(\Omega)},$$
and
$$E(u_0)\leq 0 \quad \mbox{and}\quad u_0\not\equiv 0$$
Then for all $t>0$ (where the solution exists), we have
$$u(t)> -||u(t)||_{L^2(\Omega)}.$$
\end{lem}

\noindent\textbf{Proof of Lemma \ref{pro:4}}\\
Let us define the set for every $T>0$
$$\Sigma_T=\left\{(x,t)\in \Omega\times (0,T),\quad u(x,t)<
  -||u(t)||_{L^2(\Omega)}\right\},$$
and the function
$$v(x,t)=-||u(t)||_{L^2(\Omega)}.$$
If $\Sigma_T\not=\emptyset$, then the functions $u$ and $v$
satisfy (using the condition $p\leq 2$)
$$\left\{\begin{array}{l}
\Delta u -u_t = ||u||^p_{L^p(\Omega)}-|u|^p
<||u||^p_{L^p(\Omega)}-||u||^p_{L^2(\Omega)}\leq 0\\
\\
\Delta v -v_t = \frac{F'(t)}{2\sqrt{F(t)}}\geq 0
\end{array}\right| \quad \mbox{on}\quad \Sigma_T$$
where we have used the fact that $F'\geq 0$ if $E(u_0)\leq 0$ (see
Lemma \ref{pro:2}). Consequently we have for $w=u-v$:
$$\left\{\begin{array}{l}
\Delta w -w_t <0   \quad \mbox{on}\quad \Sigma_T   \\
\\
w=0  \quad \mbox{on}\quad \left(\partial\Sigma_T\right)\backslash
\left\{t=T\right\}.
\end{array}\right.$$
The maximum principle implies $w\geq 0$ on $\Sigma_T$ which gives a
contradiction with the definition of $\Sigma_T$. Therefore
$\Sigma_T=\emptyset$ for every $T>0$, and then $w$ satisfies
$$w\geq 0 \quad \mbox{on}\quad \Omega\times (0,+\infty ).$$
From Lemma \ref{pro:2}, if $E(0)\leq 0$ and $u_0\not\equiv 0$,
then
$$\frac{F'(t)}{2\sqrt{F(t)}}\geq
\sqrt{F(0)}\frac{(p-1)\lambda_1(\Omega)}{2}e^{\frac{(p-1)\lambda_1(\Omega)}{2}
t}>0.$$ Then, if there is a point $P=(x_0,t_0)\in
\overline{\Omega}\times (0,+\infty )$ such that $w(P)=0$, we have
$\Delta w (P)-w_t(P)=-\frac{F'(t_0)}{2\sqrt{F(t_0)}}<0$, and then
there is a connected open neighborhood $\sigma_P$ of $P$ in
$\Omega\times (0,+\infty)$ such that
\begin{equation}\label{eq:12}
\left\{\begin{array}{l}
\Delta w -w_t <0  \quad \mbox{on}\quad \sigma_P   \\
\\
w\geq 0 \quad \mbox{on}\quad \sigma_P   \\
\\
w(P)=0.
\end{array}\right.
\end{equation}
As a consequence of the strong maximum principle, we get that
$w=0$ on $\sigma_P$ which does not satisfies the parabolic
equation (\ref{eq:12}). Contradiction. We conclude that
$$w>0 \quad \mbox{on}\quad \Omega \times (0,+\infty)$$
$\hfill\square$

\begin{lem}\label{pro:6}\textbf{ (Monotonicity of the infimum of $\mathbf{u}$
 for $\mathbf{p\in (1,+\infty)}$)}\\
Let us consider $p\in (1,+\infty)$. We assume that there exists $0\leq t_1< t_2$, such that
\begin{equation}\label{eq:13}
\inf_x u(t)< -||u(t)||_{L^p(\Omega)} \quad \mbox{for all}\quad t\in (t_1,t_2)
\end{equation}
and $u\in L^\infty_{loc}((0,t_2);L^p(\Omega))$.
Then the infimum
$$m(t)=\inf_x u(t)$$
is nondecreasing on $(t_1,t_2)$.
\end{lem}

\noindent\textbf{Proof of Lemma \ref{pro:6}}\\
For every $t_0\in (t_1,t_2)$, let us consider the solution $g^{t_0}=g$ of
the following ODE:
$$\left\{\begin{array}{l}
g'(t)=|g|^p-||u(t)||_{L^p(\Omega)}^p \quad \mbox{on}\quad (t_0,t_2)\\
\\
g(t_0)=m(t_0),
\end{array}\right.$$
and the set
$$\Sigma=\left\{(x,t)\in \Omega\times (t_0,t_2),\quad
  u(x,t)<g^{t_0}(t)\right\}$$
If $\Sigma\not=\emptyset$, then $u$ satisfies
$$\left\{\begin{array}{l}
\displaystyle{\Delta u-u_t\leq ||u(t)||_{L^p}^p-|g^{t_0}|^p =\Delta g^{t_0}
-\left(g^{t_0}\right)_t  \quad \mbox{on}\quad  \Sigma}\\
\\
u=g^{t_0} \quad \mbox{on}\quad \partial\Sigma.
\end{array}\right.$$
Therefore the maximum principle implies that $u\geq g^{t_0}$ on $\Sigma$,
which gives a contradiction with the definition of $\Sigma$. Thus
$\Sigma\not=\emptyset$ and $u\geq g^{t_0}$ on $\Omega\times (t_0,t_2)$,
which implies
$$m(t)\geq g^{t_0}(t) \quad \mbox{for all}\quad t\in (t_0,t_2).$$
Now using (\ref{eq:13}), we get
\begin{equation}\label{eq:16}
(g^{t_0})'(t)=|g^{t_0}(t)|^p-||u(t)||_{L^p(\Omega)}^p\geq
|m(t)|^p-||u(t)||_{L^p(\Omega)}^p\geq 0
\end{equation}
and then for $t_0'$ satisfying $t_1<t_0<t_0'< t_2$, we get
\begin{equation}\label{eq:17}
m(t_0')\geq g^{t_0}(t_0')\geq g^{t_0}(t_0)=m(t_0)
\end{equation}
$\hfill\square$

\noindent \textbf{Proof of Theorem \ref{th:5}}\\
We assume that $E(t)\geq -C_0$ on $(0,T_0)$. Then we compute
$$\begin{array}{ll}
\displaystyle{\frac12 F'(t)} & =\displaystyle{\int_\Omega uu_t}\\
\\
& \leq \displaystyle{\frac12\left(\int_\Omega u^2 +u_t^2\right)}\\
\\
&= \displaystyle{\frac12\left( F(t) -E'(t)\right)}
\end{array}$$
We deduce
$$\begin{array}{ll}
(F+E+C_0)'(t) & \leq F(t)\\
\\
&  \leq (F+E)(t) +C_0
\end{array}$$
Consequently for $t\in (0,T)$ we get
$$F(t)\leq (F+E)(t)+C_0\leq (F(0)+E(0)+C_0)e^t$$
which proves that $F(t)$ is bounded on $[0,T)$. In particular $F$
can not blow up at time $T$. ${}\hfill\square$

\noindent\textbf{Proof of Theorem \ref{th:4}}\\
We assume that $E(u_0)\leq 0$ and $u_0\not\equiv 0$.\\
\textbf{First case: $\mathbf{u_0\geq -||u_0||_{L^2(\Omega)}}$}\\
Then by lemma \ref{pro:4}, we get $u(t)\geq
-||u(t)||_{L^2(\Omega)}$ for all $t>0$. For some $\beta\in (2,3)$,
let us consider the real $\lambda_\beta$ given by lemma
\ref{pro:3}. Then by lemma \ref{pro:2}, there exists a time
$t_\beta\geq 0$, such that $||u(t)||_{L^2(\Omega)}\geq
\lambda_\beta>0$ for every $t\ge t_\beta$. Using lemma \ref{pro:3}
(and the monotonicity of the energy $E$ given by proposition
\ref{pro:1}), we get for $t\geq t_\beta$
$$F'(t)\geq 2C_\beta F^{\frac{p+3}{4}}(t)$$
which blows up in finite time $T>0$.\\
\textbf{Second case: $\mathbf{\inf_x u_0< -||u_0||_{L^2(\Omega)}}$}\\
We know by lemma \ref{pro:6} that $m(t)=\inf_x u(t)$ is
nondecreasing as long as
$$\inf_x u(t)< -||u(t)||_{L^2(\Omega)}\leq -||u(t)||_{L^p(\Omega)}$$
because we have
$$p\leq 2$$
Then
$$m(0)\leq m(t)\leq -||u(t)||_{L^2(\Omega)}^2$$
Lemma \ref{pro:2} proves that there is necessarily one time $t_0$
such that $\inf_x u(t_0) = -||u(t_0)||_{L^2(\Omega)}^2$. We can
then apply the first case with initial time $t_0$. $\hfill\square$

\bigskip
Let us conclude this section with a partial result in the case $p>2$:\\
\begin{pro}\label{pro:5}\textbf{($\mathbf{L^\infty}$ bound from
below for $\mathbf{p\in (1,+\infty)}$)}\\
Let $p\in (1,+\infty)$, and let $u$ be a solution of the problem (\ref{eq:1})-(\ref{eq:2})
with the initial data $u(t=0)=u_0$ satisfying
$$u_0\ge -||u_0||_{L^p(\Omega)},$$
and
$$E(u_0)\leq 0 \quad \mbox{and}\quad u_0\not\equiv 0.$$
Then for all $t>0$ (where the solution exists), we have
$$u(t)> -\sup_{s\in [0,t]}||u(s)||_{L^p(\Omega)}.$$
\end{pro}
\noindent\textbf{Proof of Proposition \ref{pro:5}}\\
The proof is similar to the proof of proposition \ref{pro:4}, where we use
the function $v(x,t)=-\sup_{s\in [0,t]}||u(s)||_{L^p(\Omega)}$, which
satisfies
$$\Delta v -v_t\geq 0\quad \mbox{on}\quad \Omega \times (0,+\infty)$$
In the last part of the proof, we remark that $w=0$ on $\sigma_P$,
and by connexity of $\Omega \times (0,+\infty)$, we get $w=0$ on
$\Omega \times (0,+\infty)$. The equation on $u$ implies
$u(x,t)=constant$ on $\Omega \times (0,+\infty)$ which is in
contradiction with the fact that $\int_\Omega u(t)=0$ and
$u_0\not\equiv 0$. $\hfill\square$

\section{Optimization by a gamma-convergence technique:
proof of Theorem \ref{th:optimization}}

In this whole section we assume that $\Omega$ is a bounded domain.\\



To do the proof of Theorem \ref{th:optimization}, we first need to
rewrite an integral as follows:

\begin{lem}\label{lem:10}\textbf{(Rewrite $\mathbf{\int_\Omega v|v|^p}$
as the integral of nonnegative function, for $\mathbf{p>1}$)}\\
Let us denote
$$\displaystyle A:=\left\{v\in L^{p+1}(\Omega);\quad \int_\Omega v=0;\quad \int_\Omega
v^2=1;\quad v\ge -1\ \mbox{on}\ \Omega\right\}.$$
Then there
exists a function $f\in C^2([-1,+\infty))$ which satisfies
$$f>0\quad \mbox{on}\quad (-1,1)\cup (1,+\infty),\quad \mbox{and}\quad f(-1)=f(1)=0$$
such that for every $v\in A$ we have
$$\int_\Omega v|v|^p= \int_\Omega f(v) \geq 0.$$
\end{lem}

\noindent\textbf{Proof of Lemma \ref{lem:10}}\\
Here we use the function
$$f(v)=v|v|^p-v +\frac{p}{2}(1-v^2)$$
and use the fact that $\int_\Omega v=\int_\Omega (1-v^2)=0$.
The properties of this function can be easily checked, computing
$$f'(v)=(p+1)|v|^p-pv-1,\quad f''(v)=p(p+1)v|v|^{p-2}-p.$$
${}\hfill\square$

\noindent\textbf{Proof of Theorem \ref{th:optimization}}\\
To prove Theorem \ref{th:optimization}, we simply observe that for
$\theta =\varepsilon^2$, and $v\in A$, we can write
$$\int_\Omega v|v|^p+\theta |\nabla v|^2 = \varepsilon J^\varepsilon(v)$$
with
$$J^\varepsilon(v)=\left\{\begin{array}{ll}
\displaystyle{\int_\Omega \varepsilon |\nabla v|^2 +\frac{1}{\varepsilon} f(v)} &\quad
\mbox{if}\quad v \in A\cap H^1(\Omega)\\
\\
+\infty &\quad \mbox{if}\quad v \not\in A\cap H^1(\Omega).\\
\end{array}\right.$$
As $\varepsilon$ goes to zero, the minimizers of $J^\varepsilon$ will
concentrate on the minima of the function $f$, namely on the values
$v=-1$ or $v=1$. We see formally that at the limit, we will
get discontinuous functions. To perform rigorously the analysis, we need to
introduce the space $BV(\Omega)$ of functions of bounded variations on
$\Omega$.\\
For a function $v\in L^1(\Omega)$, we define the total variation of $v$ as
$$|\nabla v|(\Omega)=\sup\left\{\int_\Omega -\sum_{i=1}^n v \frac{\partial
    \phi_i}{\partial x_i},\quad \phi=(\phi_1,...,\phi_n)\in
    (C^1(\overline{\Omega}))^n,\quad \sum_{i=1}^n \phi_i^2
\leq 1 \quad\mbox{on}\quad \Omega\right\}.$$ Then the norm in
$BV(\Omega)$ is defined by
$$\|v\|_{BV(\Omega)}:=\int_\Omega|v|\,dx+|\nabla v|(\Omega)$$
and the space $BV(\Omega)$ is naturally defined by
$$BV(\Omega)=\left\{v\in L^1(\Omega),\quad \|v\|_{BV(\Omega)}<+\infty
\right\}$$ It is known that $BV(\Omega)$ is a Banach space.

Then Theorem \ref{th:optimization} is a consequence of the
following result:
\begin{pro}\label{pro:10bis}\textbf{ (Limit inf of the energy)}\\
Assume that $\Omega$ is a bounded domain and  $\varepsilon >0$. Then the energy
$$I^\varepsilon=\inf_{u\in A} J^\varepsilon(u)$$
satisfies
$$\liminf_{\varepsilon\to 0}I^\varepsilon \geq  I^0>0$$
where $I^0$ is a constant.
\end{pro}

We give the sketch of the proof of proposition \ref{pro:10bis}
below, based on a Gamma-convergence technique, but for the
convenience of the reader, we provide in the appendix  a
self-contained proof (see Proposition \ref{pro:10} and its proof).

\noindent \textbf{Proof of Proposition \ref{pro:10bis}}\\
We remark that
$$\inf_{u\in A} J^\varepsilon(u)\geq \inf_{u\in A_0} J^\varepsilon(u)$$
where
$$A\subset \displaystyle A_0:=\left\{v\in L^{1}(\Omega);\quad \int_\Omega
  v=0;\quad v\ge -2\ \mbox{on}\ \Omega\right\}$$
with the function $f$ extended on $[-2,-1]$ by $f(v)=|v+1|$. Now
we apply the result of Modica \cite[Theorem I page 132,
Proposition 3 page 138]{modica}, with $W(t)=f(t-2)$,
$\alpha=1,\beta =3,m=2, |\Omega|=1$, $k=p+1$. It is easy to see
that there exists a constant $I_0>0$ as stated in Proposition
\ref{pro:10bis}. ${}\hfill\square$

See also the overview of Alberti \cite{Alb},
where the full Gamma-convergence result is stated.
The concept of Gamma-convergence has been introduced by De Giorgi
\cite{DeG}, and one of the first illustration of this concept was the work
of Modica, Mortola \cite{MM1}. For an introduction to Gamma-convergence and
many references, we refer the reader to the book of Dal Maso \cite{Dal}.


\section{Explicit global existence for $p=2$
and proof of Theorem \ref{th:global}}

In this section, in order to make clear the dependence on the volume
$|\Omega|$, we do not assume $|\Omega|=1$.\\

\noindent Theorem \ref{th:global} is actually a special case of
the following more general result:

\begin{theo}\label{th:globalbis}\textbf{ (Global existence for small initial
data with explicit constants, case $p=2$)}. Let $r$ be any real
number satisfying $r>\frac N2$ and $r\ge 2$. Let $\rho_r(\Omega)$
be the constant given by
$$\rho_r(\Omega)=\frac{\lambda_1(\Omega)}{4r}\cdot
\exp\left(-\left(\gamma_{N,r}\, \left(\lambda_1(\Omega)\right)^{\frac{N}{2}}\,
      |\Omega|\,\, H(\Omega)\right)^{\frac{2}{N}}\right),$$
where
$$\gamma_{N,r}=\left(1+\frac{e^{-1}}{2(1-\frac{N}{2r})}\right)^r \, 2^{r-1}\,
r^{-\frac{N}{2}}$$
For every $u_0\in
C(\overline{\Omega})$ satisfying $\int_\Omega u_0=0$ and
\begin{equation}\label{hypu0s}
\|u_0\|_{L^\infty(\Omega)}\le  \rho_r(\Omega),
\end{equation}
the (unique) solution of the problem (\ref{eq:1})-(\ref{eq:2}) is
defined for all $t$. Moreover for all $t\ge 1$ the solution
satisfies:
$$\|u(t)\|_{L^\infty(\Omega)}\leq C(r)\|u_0\|_{L^\infty}\,
e^{-\frac {\lambda_1(\Omega)}{r}t},$$ where
$$C(r):=2^{r-1\over r}|\Omega|^{\frac 1r} H(\Omega)^{\frac 1r}
\left[ 1+\frac {2\|u_0\|_{L^\infty}}{1-\frac N{2r}} \right].$$
\end{theo}

\bigskip

To deduce Theorem \ref{th:global}, we simply apply Theorem
\ref{th:globalbis} with $r=2N$, and use the fact that
$\left(\gamma_{N,2N}\right)^{\frac{2}{N}} \leq \gamma_N$ and the
inequality (\ref{iso}).

\bigskip

\noindent\textbf{Proof of Theorem \ref{th:globalbis}}

Let us denote by $K(t,x,y)$ the heat kernel associated to the
Laplacian in $\Omega$ with Neumann boundary conditions\footnote{In
fact for this section it suffices to consider  a domain $\Omega$
satisfying the extension property (see \cite[Section
1.7]{Davies}).}. That is
$$\left\{\begin{array}{ll}\frac{\partial}{\partial
t}K(t,x,y)-\Delta_x K(t,x,y)=0,&\forall x,y\in\Omega,\,t>0\\
\frac{\partial }{\partial n_x}K(t,x,y)=0,&\forall x\in
\partial\Omega,\,y\in \Omega,\,t>0\\
K(t,x,y)\longrightarrow \delta_y(x) \; \mbox{in}\; \mathcal{
D}'(\Omega)\; \mbox{as}\;  t\to0^+, & \forall y\in \Omega .
\end{array}\right.$$
This function is related to the eigenvalues and eigenfunctions of
the Neumann Laplacian $-\Delta$ in  $ \Omega$ by the following
identity:
\begin{equation}\label{Kdev}
K(t,x,y)=\sum_{k\ge 0} e^{-\lambda_k (\Omega) t} f_k(x)f_k(y)
={1\over |\Omega|}+\sum_{k\ge 1} e^{-\lambda_k (\Omega) t}
f_k(x)f_k(y),
\end{equation}
where $\left\{\lambda_k(\Omega)\; ; \; k\ge 0\right\} $ are the
eigenvalues of $-\Delta$ and $\left\{f_k\; ; \; k\ge 0\right\}$ is
an $L^2$-orthonormal family of corresponding eigenfunctions
(recall that $\lambda_0(\Omega)=0$ and $f_0={1\over
|\Omega|^{1\over 2}}$). Let us set $ K_0(t,x,x) := K(t,x,x)
-{1\over |\Omega|} $. From the classical results on heat kernels
(see for instance \cite[Theorem 2.4.4]{Davies}) we know that the
function $t^{N/2}K(t,x,x)$ is bounded on $(0,1]\times\Omega$, and
then the same is true for $t^{N/2}K_0(t,x,x)$. On the other hand,
it follows immediately from (\ref {Kdev}) that
$e^{\lambda_1(\Omega) (t-1)} K_0(t,x,x) $ is decreasing on $[1,
+\infty)$ and then, for any $t\ge 1$,
$$K_0(t,x,x)  \le e^{-\lambda_1(\Omega) (t-1)} K_0(1,x,x) \le
C_1 e^{-\lambda_1(\Omega) (t-1)} \le C_2 t^{-N/2}$$
for some constants $C_1$ and $C_2$. Hence, $t^{N/2}K_0(t,x,x)$ is
bounded on $(0,+\infty)\times\Omega$ and we denote by $H(\Omega)$
its supremum. Since for all $t>0$, $K(t,x,y)$ achieves its
supremum on the diagonal of $\Omega\times\Omega$, the constant
$H(\Omega)$ is actually the best constant in the following
inequality
$$K(t,x,y)\le {1\over |\Omega|}+H(\Omega)t^{-N/2}$$
valid in $(0,+\infty)\times\Omega\times\Omega$. Notice that, in
contrast to the Dirichlet boundary condition case, there is no
universal upper bound to $H(\Omega)$ (even for domains of fixed
volume). Indeed, it is rather easy to see that, in the Dirichlet
case, the heat kernel is bounded above by $(4\pi t)^{-N/2}$
whatever the domain is.

\begin{rem}\label{rem:H}
Notice that
$$\displaystyle H(\Omega)\ge {t^{N\over 2}\over
|\Omega|}\left[\left(\int_\Omega K(t,x,x)\right) -1\right]={t^{N\over 2}\over
|\Omega|}\left[ \hbox{tr}\, e^{t\Delta} -1\right]$$
with
$$\displaystyle \hbox{tr}\, e^{t\Delta}=(4\pi t)^{-N\over
2}\left[|\Omega|+{\sqrt{\pi}\over 2}\left({\mathcal H}^{N-1}\left(\partial\Omega\right)\right)t^{1\over 2}
+o(t)\right]\quad\mbox{as}\quad t\to 0^+,$$
 see \cite{G}, see also \cite{AM},
\cite{Bu} and \cite{Me} for the dependance of $\lambda_1(\Omega)$
on the geometry of $\Omega$. This implies
$$H(\Omega)\geq (4\pi)^{-N/2}.$$
\end{rem}

\bigskip

The following property seems to be a standard one, we give it for
completeness.

\begin{lem}\textbf{ ($\mathbf{L^p-L^q}$-estimate for the linear
heat equation)}

Let $v_0\in C(\overline{\Omega})$ be such that
$\int_{\Omega}v_0=0$ and let
$$v(t,x)=\left(e^{t\Delta}v_0\right)(x)=\int_\Omega
K(t,x,y)v_0(y)\,dy$$ be the solution of the heat equation with
Neumann boundary condition and $v_0$ as initial data. For any
positive $t$ and all $1<p\le q\le +\infty$, we have
\begin{equation}\label{LpLq}
\|v(t)\|_{L^q}\le {\left[ 2^{p-1}H(\Omega)t^{-N/2}\right]}^{\frac
1p-\frac 1q}\|v_0\|_{L^p}.
\end{equation}
\end{lem}
\noindent\textbf{Proof.} Since $\int_{\Omega}v_0 =0$ we have, for
any $p>1$ and any $(t,x)\in (0,+\infty)\times\Omega$,
$$|v(t,x)|=\Big|\int_\Omega K_0(t,x,y)\, v_0(y)\,dy\Big|\le \|v_0\|_{L^p}\|K_0(t,x,\cdot)\|_{L^{p'}}$$
with $\frac 1p+\frac 1{p'}=1$. Now,
\begin{eqnarray*}
\|K_0(t,x,\cdot)\|_{L^{p'}}^{p'}&=&\int_\Omega|K_0(t,x,y)|^{p'}\,dy\\
&\le&{\left[H(\Omega)t^{-N/2}\right]}^{p'-1}\|K_0(t,x,\cdot)\|_{L^1}\\
\end{eqnarray*}
with $\|K_0(t,x,\cdot)\|_{L^1} =\int_\Omega |K(t,x,y)- {1\over
|\Omega|}|\,dy \le \int_\Omega (K(t,x,y)+{1\over |\Omega|})\,dy =
2$, since $K(t,x,y)\geq 0$ and $\int_\Omega K(t,x,y)\,dy=1$.
Hence, $\|K_0(t,x,\cdot)\|_{L^{p'}}^{p'}\le 2
{\left[H(\Omega)t^{-N/2}\right]}^{p'-1}$ and then,
$$\|v(t)\|_{L^\infty}\le 2^{1\over p'}\|v_0\|_{L^p}{\left[H(\Omega)t^{-N/2}\right]}^{\frac{p'-1}{p'}}
=\|v_0\|_{L^p}{\left[2^{p-1}H(\Omega)t^{-N/2}\right]}^{\frac{1}{p}}.$$
Now let us remark that
$$\begin{array}{ll}
||v(t)||_{L^q}^q &=\quad \displaystyle{\int_\Omega |v(t)|^{q-p}|v(t)|^p} \\
\\
&\leq \quad ||v(t)||_{L^\infty}^{q-p} ||v(t)||_{L^p}^p\\
\\
& \leq \quad \|v_0\|_{L^p}^{q-p}{\left[2^{p-1}
H(\Omega)t^{-N/2}\right]}^{\frac{q-p}{p}}||v(t)||_{L^p}^p.
\end{array}$$
We get finally the $L^p-L^q$ estimate of the lemma, using the
contraction property of the Heat equation with Neumann boundary
condition (see \cite[Section 3.1.1]{Lun}): $||v(t)||_{L^p}\leq
||v_0||_{L^p}$. ${}\hfill\square$

\bigskip

The following elementary property will also be useful

\begin{lem}\label{lInt}
Let $\alpha,\beta\geq 1$. For all $f\in L^{\alpha+\beta}(\Omega)$,
we have
$$\left(\int_{\Omega} |f|^\alpha \right)\cdot  \left(\int_{\Omega}
  |f|^\beta\right) \leq
|\Omega|\int_{\Omega} |f|^{\alpha+\beta}.$$
\end{lem}
\noindent\textbf{Proof.} Using H\"older inequality, the following
holds
\begin{eqnarray*}
\int_{\Omega} |f|^\alpha &\leq& {\left[ \int_{\Omega}
|f|^{\alpha+\beta} \right]}^{\frac{\alpha} {\alpha+\beta}}
\times |\Omega|^{\frac \beta{\alpha+\beta}}\\
\int_{\Omega} |f|^\beta &\leq& {\left[ \int_{\Omega}
|f|^{\alpha+\beta} \right]}^{\frac{\beta}{\alpha+\beta}}\times
|\Omega|^{\frac {\alpha}{\alpha+\beta}}
\end{eqnarray*}
Thus
$$
\int_{\Omega} |f|^\alpha \times \int_{\Omega} |f|^\beta\leq
|\Omega|\int_{\Omega} |f|^{\alpha+\beta}
$$
${}\hfill\square$

\bigskip

From now on, we denote by $u$ a solution of
(\ref{eq:1})-(\ref{eq:2}) for $p=2$ whose initial data $u_0$
satisfies $\int_\Omega u_0 =0$. The proof of Theorem
\ref{th:global} is based on the following a priori estimate:

\begin{lem}\label{l3}{\bf (A priori estimate)}\\
Let $r>\frac N2$ satisfying $r\geq 2$, $0<K<\frac {\lambda_1}{2r}$, and assume that
there exists a positive $T$ such that $\|u(t)\|_{L^\infty}\le K$
for all $t\in[0,T]$. Then, for any $\beta\in (0,1)$, we have
\begin{eqnarray}\label{priori1}
\|u(T)\|_{L^\infty}&\le& \|u_0\|_\infty 2^{1-{\frac
1r}}H(\Omega)^{\frac1r}|\Omega|^{\frac 1r}[\beta
T]^{-\frac N{2r}} \times\nonumber\\
&&\hskip -1cm\left[e^{-2[\frac{\lambda_1}r-K] (1-\beta)T}+
\frac{\beta e^{-1}\|u_0\|_{L^\infty}} {(1-\beta)
[\lambda_1/r-2K]\left(1-\frac N{2r}\right)} \right].
\end{eqnarray}
\end{lem}
\noindent\textbf{Proof} of the lemma. For simplicity we will write
$\lambda_1$ and $H$ for $\lambda_1(\Omega)$ and $H(\Omega)$. Let
$q\geq 2$ be an even integer. For all $t\in [0,T]$, we have
$|u|^{q+1}\le Ku^q$ and then, setting $F_q(t):=\frac 1q\int
u^q(t)\,dx\ge 0 $, $\int |u|^{q+1}\le qKF_q(t)$. Using lemma
\ref{lInt} we also have
$$\left|\left(\int u^{q-1}\right) \left(\intm u^2\right)\right|\le
\left(\int |u|^{q-1}\right)\left(\intm u^2\right) \le
\int |u|^{q+1}\le qKF_q(t).$$
Therefore,
$$\left|\int u^{q+1}-\left(\int u^{q-1}\right)\left(\intm u^2\right)\right|\le 2qKF_q(t).$$
Multiplying Equation (\ref{eq:1}) (with $p=2$) by $u^{q-1}$
and integrating over $\Omega$ we get, after integrating by parts,
\begin{equation}\label{eqEp}\int
u^{q-1}u_t+4\frac{q-1}{q^2}\int |\nabla u^{q/2}|^2\le 2qKF_q(t).
\end{equation}
This implies that $F_q'(t)-2qKF_q(t)\le 0$.
Hence, for all $t\in [0,T]$, we have $F_q(t)\le F_q(0)e^{2qKt}$,
or equivalently
$$\|u(t)\|_{L^q}\le \|u(0)\|_{L^q}e^{2Kt}.$$
Making $q\to\infty$, one can deduce that
\begin{equation}\label{eqLinfty}\|u(t)\|_{L^{\infty}}\le
\|u_0\|_{L^{\infty}}\,e^{2Kt}.
\end{equation}
On the other hand, we clearly have $\int u(t)=0$ for all $t$.
Poincar\'e's inequality gives, for all $t\in (0,T]$, $\int|\nabla
u|^2\ge \lambda_1\int u^2$. Taking $q=2$ in (\ref{eqEp}) we then
get, for all $t\in[0,T]$,
\begin{equation}\label{eqL2}\|u(t)\|_{L^2}\le
\|u_0\|_{L^2}\,e^{-(\lambda_1-2K)t}\le
\|u_0\|_{L^\infty}|\Omega|^{\frac 12}e^{-(\lambda_1-2K)t}.
\end{equation}
By $L^2-L^\infty$ interpolation, we obtain, since $2\le r< +\infty$
\begin{eqnarray}\label{eqLr}\|u(t)\|_{L^r}\le|\Omega|^{\frac 1r}\|u_0\|_{L^\infty}
e^{-2[\frac {\lambda_1}r-K]t}.
\end{eqnarray}
Now, Equation (\ref{eq:1}) leads to the following integral
equation
$$u(t+t_0)=e^{t_0\Delta}u(t)+\int_0^{t_0}e^{(t_0-s)\Delta}f(u(t+s))\,ds,$$
with $f(u):=u^2-\intms u^2$. Taking $t_0:=\beta T$ and
$t=t_1:=(1-\beta)T$ and using the $L^r-L^\infty$ estimates
(\ref{LpLq}) we get
\begin{eqnarray*}
\|u(T)\|_{L^\infty}&=&\|u(t_1+t_0)\|_{L^\infty}\le\|e^{t_0\Delta}u(t_1)\|_{L^\infty}+
\int_0^{t_0} \|e^{(t_0-s)\Delta}f(u(t_1+s))\|_{L^\infty}\,ds\\
&\le& 2^{1-{\frac 1r}}H^{\frac
1r}t_0^{-\frac{N}{2r}}\|u(t_1)\|_{L^r}+ \int_0^{t_0} 2^{1-{\frac
1r}}H^{\frac 1r}(t_0-s)^{-\frac{N}{2r}}\|f(u(t_1+s))\|_{L^r}\,ds.
\end{eqnarray*}
Using H\"{o}lder inequality for the first line and
(\ref{eqLr}) for the second, we get
\begin{equation}\label{eq:2.15b}
\begin{array}{lll}\|f(u(t))\|_{L^r}&\le&
2\|u(t)\|_{L^{2r}}^2\\
&\le&2|\Omega|^{\frac 1r} \|u_0\|_{L^\infty}^2 e^{-2
[\frac{\lambda_1}r-2K]t}.
\end{array}
\end{equation}
Setting $\alpha_1:=2[\frac{\lambda_1}r-K]> 0$ and
$\alpha_2:=2[\frac{\lambda_1}r-2K]> 0$ we obtain using (\ref{eqLr}) and
(\ref{eq:2.15b}):
\begin{eqnarray*}
\frac{\|u(T)\|_{L^\infty}}{\left(2^{1-{\frac 1r}}H^{\frac 1r}\right)} &
\le & t_0^{-\frac{N}{2r}}\|u(t_1)\|_{L^r}+
\int_0^{t_0}(t_0-s)^{-\frac{N}{2r}} \|f(u(t_1+s))\|_{L^{r}}\,ds \\
&\le& |\Omega|^{\frac 1r}\|u_0\|_{L^\infty} t_0^{-\frac{N}{2r}}e^{-\alpha_1
  t_1} \\
&& +2|\Omega|^{\frac 1r}\|u_0\|_{L^\infty}^2 \int_0^{t_0}  (t_0-s)^{-\frac
N{2r}} e^{-\alpha_2(t_1+s)}\,ds \\
&\le & |\Omega|^{\frac 1r}\|u_0\|_{L^\infty}
t_0^{-\frac{N}{2r}}e^{-\alpha_1 t_1}\\
&&+2|\Omega|^{\frac 1r} \|u_0\|_{L^\infty}^2e^{-\alpha_2t_1} \int_0^{t_0}
(t_0-s)^{-\frac N{2r}} \,ds\\
&=& |\Omega|^{\frac 1r} \|u_0\|_{L^\infty} t_0^{-\frac{N}{2r}}
\left[e^{-\alpha_1t_1}+2\|u_0\|_{L^\infty}
\frac{t_0 e^{-\alpha_2t_1}}{1-\frac N{2r}} \right]\\
&\le& |\Omega|^{\frac 1r}\|u_0\|_{L^\infty} t_0^{-\frac{N}{2r}}
\left[e^{-\alpha_1(1-\beta)T}+
\frac{2\|u_0\|_{L^\infty}\beta T e^{-\alpha_2(1-\beta)T}
}{1-\frac N{2r}} \right]\\
&\le& |\Omega|^{\frac 1r} \|u_0\|_{L^\infty} t_0^{-\frac{N}{2r}}
\left[e^{-\alpha_1(1-\beta)T}+
\frac{2\beta\|u_0\|_{L^\infty}}{1-\frac N{2r}}
\sup_{\overline{T}>0}\left(\overline{T}
  e^{-\alpha_2(1-\beta)\overline{T}}\right)\right]\\
&\le& |\Omega|^{\frac 1r}\|u_0\|_{L^\infty}(\beta
T)^{-\frac{N}{2r}}
\left[e^{-\alpha_1(1-\beta)T}+\frac{2\beta
e^{-1}\|u_0\|_{L^\infty}}{1-\frac N{2r}}\cdot\frac
1{\alpha_2(1-\beta)} \right].
\end{eqnarray*}
${}\hfill\square$

\noindent\textbf{End of the proof of Theorem \ref{th:globalbis}}\\
\underline{First step}: Global existence.\\
Let $u$ be a solution of (\ref{eq:1})-(\ref{eq:2}) whose initial
data $u_0$ satisfies (\ref{hypu0s}). Let us suppose, for a
contradiction, that the maximal time of existence $T_{max}$ of $u$
is finite. Put $K=\frac {\lambda_1}{4r}$ in
the last lemma. Let $T$ be the maximal time such that
$\|u(t)\|_{L^\infty}\le K$ in $[0,T]$. Hence (see
(\ref{eqLinfty}))
\begin{equation}\label{star}
K=\|u(T)\|_{L^\infty}\le \|u_0\|_{L^\infty}e^{2KT}.
\end{equation}
From the assumption (\ref{hypu0s}), we have $\|u_0\|_{L^\infty}
\le \alpha \frac {\lambda_1}{4r}=\alpha K$, with
$$\alpha=\exp\left(-\left(\gamma_{N,r}\,  H|\, \Omega| \, \lambda_1^{\frac{N}{2}}\right)^{\frac{2}{N}}\right) \quad < \quad 1$$
After
replacing into (\ref{star}), this gives
\begin{equation}\label{Tlower}
T\ge-\frac{\ln\alpha}{2K}.\end{equation} Applying Lemma \ref{l3} with $\beta=\frac 12$,
we get:
$$
\|u(T)\|_{L^\infty}\le 2^{r-1\over r}H^{\frac
1r} |\Omega|^{\frac 1r}(T/ 2)^{-\frac{N}{2r}}\|u_0\|_{L^\infty}\left[e^{-\frac {3
\lambda_1}{4r}T}+\frac{ 2 e^{-1}\|u_0\|_{L^\infty}}{\frac
{\lambda_1}{r} \left(1-\frac N{2r}\right)} \right].$$
Since
$\|u_0\|_{L^\infty} \le \alpha \frac {\lambda_1}{4r} =\alpha K$, it follows
that
\begin{eqnarray*}
\frac{\|u(T)\|_{L^\infty}}{K}&\le &
2^{r-1\over r}H^{\frac 1r}|\Omega|^{\frac 1r}(T/2)^{-\frac{N}{2r}} \alpha
\left[e^{-\frac {3 \lambda_1}{4r}T}+\frac{e^{-1}\alpha}{ 2\left(1-\frac N{2r}\right)} \right]\\
&<&2^{r-1\over r}H^{\frac 1r}|\Omega|^{\frac 1r}
(T/2)^{-\frac{N}{2r}}
\left[1+\frac{e^{-1}}{ 2\left(1-\frac N{2r}\right)} \right]\\
&\le&2^{r-1\over r}H^{\frac 1r}|\Omega|^{\frac 1r} \left(\frac{(-\ln
    \alpha)}{\frac{\lambda_1}{r}}\right)^{-\frac{N}{2r}}\left[1+\frac{e^{-1}}{ 2\left(1-\frac N{2r}\right)} \right]\\
&=&\left(\gamma_{N,r}H |\Omega| \lambda_1^{\frac{N}{2}} (-\ln
  \alpha)^{-\frac{N}{2}}\right)^{\frac1r} \\
&=& 1
\end{eqnarray*}
which shows that
$$\|u(T)\|_{L^\infty} < K.$$
This contradicts the definition of $T$.\\
\underline{Second step}: Exponential decay.\\
Note that, since $K= \frac{\lambda_1}{4r} $, estimate
(\ref{eqLr}) becomes
$$\|u(t)\|_{L^{r}}\le |\Omega|^{\frac 1r} \|u_0\|_{L^\infty}\,
e^{-\frac{3\lambda_1}{2r}t} $$ for all $t\in[0,+\infty[$. Using
again the integral equation, we have
$$u(t+1)=e^{\Delta}u(t)+\int_0^{1}e^{(1-s)\Delta}f(u(t+s))\,ds.$$
Using the
same computation as in the proof of Lemma \ref{l3}, we obtain for all $t>0$
\begin{eqnarray*}
\left(2^{r-1\over r}H^{\frac 1r}\right)^{-1}\|u(t+1)\|_{L^\infty} &\le &
\|u(t)\|_{L^r}
+\int_0^1 (1-s)^{-\frac N{2r}}\|f(u(t+s))\|_{L^r}\,ds\\
&\le& |\Omega|^{\frac 1r}\|u_0\|_{L^\infty}e^{-\frac{\lambda_1}{r}t}
\left[ 1+2\|u_0\|_{L^\infty}\int_0^1(1-s)^{-\frac N{2r}}\right]\,ds\\
&\le& \left(2^{r-1\over r}H^{\frac 1r}\right)^{-1}C(r)\|u_0\|_{L^\infty}e^{-\frac{\lambda_1}{r}t}.
\end{eqnarray*}
${}\hfill\square$


\section{Appendix : proof of a Gamma-convergence-like result}

We give here the following result which is more precise than
proposition \ref{pro:10bis}, and propose a self-contained proof.

\begin{pro}\label{pro:10}\textbf{ (Limits of the energy of the minimizers)}\\
Assume that $\Omega$ is a bounded domain. Then for every
$\varepsilon >0$, there exists at least one minimizer
$u^\varepsilon$ of the following problem
$$I^\varepsilon=\inf_{u\in A} J^\varepsilon(u)$$
and
$$\liminf_{\varepsilon\to 0}I^\varepsilon \geq J^0(u^0)\geq I^0>0$$
More precisely, there exists a subsequence
$(u^{\varepsilon'})_{\varepsilon'}$ such that
$$u^{\varepsilon'}\longrightarrow u^0 \quad \mbox{in}\quad L^1(\Omega)$$
and $u^0\in B$, where
$$B:=\left\{u\in BV(\Omega);\quad u=\pm 1\mbox{ a.e. in }\Omega;
\quad \int_\Omega u=0 \right\}$$
and 
$$I^0=\inf_{u\in B}J^0(u)$$
where
$$J^0(u)=c|\nabla u|(\Omega) \quad \mbox{with}\quad c=\int_{-1}^1 \sqrt{f(s)}\, ds.$$
\end{pro}

To prove proposition \ref{pro:10}, we will use the following
classical compactness result in $BV(\Omega)$.

\begin{pro}\label{pro:11}\textbf{ (Compactness in $BV(\Omega)$, \cite{giusti})}\\
Let $\Omega$ be a bounded domain. For every sequence $(v^n)_n$,
bounded in $BV(\Omega)$, there exists a subsequence
$(v^{n'})_{n'}$ and $v^\infty\in BV(\Omega)$ such that
$$v^{n'}\longrightarrow v^\infty\quad\mbox{in}\quad L^1(\Omega)$$
and
$$\liminf_{n'\to +\infty }|\nabla v^{n'}|(\Omega) \quad \ge \quad |\nabla v^\infty|(\Omega).$$
\end{pro}


\noindent\textbf{Proof of Proposition \ref{pro:10}}\\
The proof of this proposition is done in the following steps.\\

\noindent \textbf{Step 1: there exists a minimizer $u^\varepsilon$.}\\
It is easy to see that there exists a constant $C>0$ such that
\begin{equation}\label{eq:20}
|v|^{p+1} +C \geq f(v)\geq |v|^{p+1} -C\quad \mbox{for}\quad v\geq
-1.
\end{equation}
Therefore, for $\varepsilon>0$ fixed, every minimizing sequence of
$J^\varepsilon$ in $A$ is bounded in $H^1(\Omega)\cap
L^{p+1}(\Omega)$. From the compactness of the injection
$H^1(\Omega)\longrightarrow L^2(\Omega)$, it is classical to get
the existence of a minimizer $u^\varepsilon\in A$ of
$J^\varepsilon$.\\

\noindent \textbf{Step 2: there exists $C_0>0$ such that
$J^\varepsilon(u^\varepsilon)\le C_0$ for $\varepsilon$ small enough.}\\
Here for $\varepsilon$ small enough we will build a function
$w^\varepsilon\in A$ such that $J^\varepsilon(w^\varepsilon)\leq C_0$.\\
Let us consider the direction $x_1$ and assume that the hyperplane
$\left\{x_1=0\right\}$ separates $\Omega$ in two equal volumes:
$$\left|\Omega\cap\left\{x_1< 0\right\}\right|=\left|\Omega\cap\left\{x_1>
    0\right\}\right|.$$
For $\delta >0$, we define the function
$$v^\delta(x_1)=\left\{\begin{array}{ll}
-1 &\quad\mbox{if}\quad x_1<-\varepsilon\\
\\
\displaystyle{\frac{x_1}{\varepsilon}} &\quad\mbox{if}\quad
-\varepsilon\leq x_1 \leq \varepsilon (1+\delta)\\
\\
1+\delta &\quad\mbox{if}\quad x_1> \varepsilon (1+\delta).
\end{array}\right.$$
Next we define the translation of $v^\delta$, for a small
parameter $a \in \R$:
$$v^\delta_a(x_1)=v^\delta(x_1-a).$$
Then for $a$ close enough to zero and fixed, there is a unique
$\delta=\delta(a)>0$ such that
$$\int_\Omega |v^{\delta(a)}_a|^2=1.$$
In particular, this implies that
$$\int_{\Omega\cap\left\{x_1>a+\varepsilon\right\}}\left(
|v^{\delta(a)}_a|^2-1\right)=
\int_{\Omega\cap\left\{a-\varepsilon<x_1<a+\varepsilon\right\}}
\left(1-|v^{\delta(a)}_a|^2\right)$$ and then
$$\frac{1}{2}\left|\Omega\right|(1- o(|a|+\varepsilon
(1+\delta)))\left((1+\delta)^2-1\right) \leq 2\varepsilon
(\mbox{diam}(\Omega))^{n-1}$$ i.e. for $a$ close enough to zero,
there exists a constant $C>0$ such that
$$\delta(a)\leq C\varepsilon \leq 1 \quad \mbox{for}\quad \varepsilon
\quad\mbox{small enough}.$$ Then we consider the map
$$a\longmapsto \Phi(a)=\int_\Omega v^{\delta(a)}_a.$$
We have  $\Phi(-2\varepsilon)>0$. On the other hand, because the
open set $\Omega$ is connected, there exists a constant $C_2>0$
such that $|\Omega\cap\left\{0<x_1<\eta\right\}|\geq C_2 \eta$ for
$\eta>0$ small enough. Therefore $\Phi(a)\leq |\Omega|/2 \delta(a)
- C_2(a-\varepsilon)$ and then $\Phi(\varepsilon
(1+C|\Omega|/C_2))<0$. Using the continuity of the map $\Phi$, we
deduce that
$$\exists a\in (-2\varepsilon, (1+C|\Omega|/C_2)\varepsilon) \quad \int_\Omega
v^{\delta(a)}_a=0.$$ We set $w^\varepsilon=v^{\delta(a)}_a\in A$
and estimate $J^\varepsilon(w^\varepsilon)$ as follows
$$\int_\Omega \varepsilon |\nabla w^\varepsilon|^2\leq
(2+\delta(a))(\mbox{diam}(\Omega))^{n-1}$$
$$\int_\Omega \frac{1}{\varepsilon}f(w^\varepsilon)\leq  \left(\sup_{[-1,1 ]}f\right)
2(\mbox{diam}(\Omega))^{n-1} +
|\Omega|\frac{1}{\varepsilon}\sup_{[1,1+\delta(a)]}f.$$ We then
remark that
$$\sup_{[1,1+\delta]}f\leq \frac{1}{2}f''(1)(\delta)^2 +o(\delta^2)$$
because $f(1)=f'(1)=0$. We deduce that
$$\frac{1}{\varepsilon}\sup_{[1,1+\delta(a)]}f\leq \frac{1}{2}f''(1) C^2
\varepsilon + o(\varepsilon)$$ Putting all together we get the
existence of a constant $C_0>0$ such that for $\varepsilon$ small
enough we get
$$J^\varepsilon(w^\varepsilon)\leq C_0.$$
Because $w^\varepsilon\in A$, and $u^\varepsilon$ is a minimizer
of $J^\varepsilon$ on $A$, we deduce that
$$J^\varepsilon(u^\varepsilon)\leq J^\varepsilon(w^\varepsilon)\leq C_0.$$
This ends the proof of Step 2.\\

\noindent \textbf{Step 3: there exists $C_1>0$ such that
$\mathbf{||v^\varepsilon ||_{BV(\Omega)}\leq C_1}$ for
$\mathbf{\varepsilon}$ small enough.}\\
We define
$$G(s)=\left\{\begin{array}{ll}
\displaystyle{\int_{-1}^s \sqrt{f}(t) \ dt} & \quad \mbox{if}\quad s\geq -1\\
\\
0  & \quad \mbox{if}\quad s < -1
\end{array}\right.$$
and
$$v^\varepsilon=G(u^\varepsilon).$$
From (\ref{eq:20}), we get for $t\geq -1$
$$f(t)\leq |t|^{p+1} +C$$
and then there exists some constants $C,C'>0$ such that (using
$p>1$)
\begin{equation}\label{eq:21}
G(s)\leq C\left(|s|^{\frac{p+3}{2}} + |s| +1\right)\leq
C'\left(|s|^{p+1}
  + |s| +1\right)
\end{equation}
To estimate $||v^\varepsilon||_{L^1(\Omega)}$, we will first
estimate
  $||u^\varepsilon||_{L^1(\Omega)}$ and $||u^\varepsilon||_{L^{p+1}(\Omega)}$.
Because $\int_\Omega u^\varepsilon =0$, we remark that
$\int_\Omega (u^\varepsilon)^+ =\int_\Omega (u^\varepsilon)^-\leq
|\Omega|$ and then
$$||u^\varepsilon||_{L^1(\Omega)}\leq  2|\Omega|.$$
From Step 2, we have $J^\varepsilon(u^\varepsilon)\leq C_0$, and
then because of (\ref{eq:20}), we get
$$\int_\Omega |u^\varepsilon|^{p+1} - C|\Omega|\leq C_0\varepsilon$$
which gives
$$||u^\varepsilon||_{L^{p+1}(\Omega)}\leq \left(C |\Omega|+C_0\varepsilon\right)^{\frac{1}{p+1}}.$$
Putting all together we deduce from (\ref{eq:21}):
$$||v^\varepsilon||_{L^1(\Omega)}=||G(u^\varepsilon)||_{L^1(\Omega)}\leq
C'\left(||u^\varepsilon||_{L^{p+1}(\Omega)}^{p+1}+
  ||u^\varepsilon||_{L^1(\Omega)}+|\Omega|\right).$$
To estimate the whole norm in $BV(\Omega)$ we only need to
estimate $|\nabla v^\varepsilon|(\Omega)$. This is done, using the
following classical trick of Modica \cite{modica} for every $u\in
A\cap H^1(\Omega)$ (and $a^2 +b^2\geq 2ab$)
\begin{equation}\label{eq:22}
J^\varepsilon(u)=\int_\Omega \varepsilon |\nabla u^\varepsilon|^2
+\frac{1}{\varepsilon}f(u) \geq \int_\Omega 2|\nabla u|\sqrt{f(u)}
= \int_\Omega 2 |\nabla G(u)|.
\end{equation}
Applied to $v^\varepsilon$, we get
$$\int_\Omega |\nabla v^\varepsilon|=\int_\Omega |\nabla
G(u^\varepsilon)| \leq \frac12 J^\varepsilon(u^\varepsilon) \leq
\frac12 C_0$$
This proves the expected inequality and ends Step 3.\\

\noindent \textbf{Step 4: $II^0\geq 2 |\nabla v^0|(\Omega)$.}\\
Let us define
$$II^0=\liminf_{\varepsilon\to 0}J^\varepsilon(u^\varepsilon).$$
Then we extract a subsequence $(u^{\varepsilon'})_{\varepsilon'}$
such that
$$J^{\varepsilon'}(u^{\varepsilon'})\longrightarrow II^0.$$
From (\ref{eq:22}), we have
$$J^{\varepsilon'}(u^{\varepsilon'})\geq \int_\Omega 2 |\nabla
G(u^{\varepsilon'})| =\int_\Omega 2 |\nabla v^{\varepsilon'}| = 2
|\nabla v^{\varepsilon'}|(\Omega).$$ From Step 3 and the
compactness result in $BV$ (proposition \ref{pro:11}), up to
extract a new subsequence, we can assume that there exists $v^0\in
BV(\Omega)$ such that
\begin{equation}\label{eq:1v1}
v^{\varepsilon'} \longrightarrow v^0 \quad\mbox{in}\quad
L^1(\Omega)
\end{equation}
and
$$\liminf_{\varepsilon'\to 0}|\nabla v^{\varepsilon'}|(\Omega) \quad \ge
\quad |\nabla v^0|(\Omega)$$ so that
$$II^0\geq 2 |\nabla v^0|(\Omega).$$
Moreover from (\ref{eq:1v1}), and the converse Lebesgue theorem,
up to extraction of a subsequence, we can assume that
\begin{equation}\label{eq:1v2}
v^{\varepsilon'}(x) \longrightarrow v^0(x) \quad \mbox{for
a.e.}\quad x\in \Omega.
\end{equation}

\noindent\textbf{Step 5. $\mathbf{u^{\varepsilon''}\longrightarrow
u^0}$ in $\mathbf{L^{p+1}(\Omega)}$ and $\mathbf{u^0=\pm 1}$ a.e.
in
$\mathbf{\Omega}$.}\\
From Step 2, we have
$$\int_\Omega f(u^\varepsilon)\leq \varepsilon
J^\varepsilon(u^\varepsilon)\leq C\varepsilon$$ and then
$$f(u^\varepsilon)\longrightarrow 0 \quad\mbox{in}\quad L^1(\Omega).$$
Then from the converse Lebesgue theorem, there exists a function
$h\in L^1(\Omega)$, that we can always choose satisfying $h\geq
1$, such that there exists a subsequence
$(u^{\varepsilon''})_{\varepsilon''}$ with
$$f(u^{\varepsilon''}(x))\leq h(x) \quad \mbox{for a.e.}\quad x\in \Omega$$
and
\begin{equation}\label{eq:25}
f(u^{\varepsilon''}(x))\longrightarrow 0 \quad \mbox{for
a.e.}\quad x\in \Omega.
\end{equation}
Our goal is now to prove that there exists a subsequence of
$(u^{\varepsilon''})_{\varepsilon''}$ which is convergent to some
$u^0$ in
$L^{p+1}(\Omega)$, and $u^0=\pm 1$ a.e. in $\Omega$.\\
We remark that from (\ref{eq:20}), we have for $v\geq 1$,
$f^{-1}(|v|^{p+1}-C)\leq v$, and then setting $h=|v|^{p+1}-C$, we
get
$$f^{-1}(h)\leq |h+C|^{\frac{1}{p+1}}.$$
This proves that $f^{-1}(h)\in L^{p+1}(\Omega)$, and
\begin{equation}\label{eq:1v3}
|u^{\varepsilon''}|\leq 1+f^{-1}(h)\in L^{p+1}(\Omega).
\end{equation}
Now from (\ref{eq:1v2}) and the continuity of $G^{-1}$ on
$(0,+\infty )$, we have
$$u^{\varepsilon''} =G^{-1}(v^{\varepsilon''})\longrightarrow G^{-1}(v^0)=:
u^0 \quad \mbox{for a.e.}\quad x\in \Omega.$$ Moreover from
(\ref{eq:1v3}), we deduce that $u^0\in L^{p+1}(\Omega)$ and
\begin{equation}\label{eq:26}
u^{\varepsilon''}\longrightarrow u^0\quad\mbox{in}\quad
L^{p+1}(\Omega).
\end{equation}
Consequently from (\ref{eq:25}), we deduce
$$f(u^{\varepsilon''}(x))\longrightarrow 0=f(u^0(x)) \quad \mbox{for a.e.}\quad
x\in \Omega$$ and then
$$u^0(x)=\pm 1 \quad \mbox{for a.e.}\quad
x\in \Omega.$$
This ends the proof of Step 5.\\

\noindent\textbf{Step 6: $\mathbf{\int_\Omega u^0=0}$,
$\mathbf{u^0\in BV(\Omega)}$ and $\mathbf{|\nabla
v^0|(\Omega)=\frac{G(1)}2 |\nabla u^0|(\Omega)}$.}\\
(i) From (\ref{eq:26}), we get in particular that
$u^{\varepsilon''}\longrightarrow u^0$ in $L^1(\Omega)$, and then
$0=\int_\Omega u^\varepsilon \longrightarrow \int_\Omega u^0$
which proves that
$$\int_\Omega u^0 =0.$$
(ii) We have $v^0=G(u^0)\in BV(\Omega)$ and $u^0=\pm 1$ a.e. in
$\Omega$. Therefore $v^0=G(\pm 1)$ a.e. in $\Omega$. Moreover,
because each of these two functions only takes two values, we can
express $u^0$ as a function of $v^0$, i.e. (using $G(-1)=0$)
$$\displaystyle u^0=\frac{2}{G(1)} v^0-1.$$
Thus $u^0\in BV(\Omega)$ and $|\nabla u^0|(\Omega)= \frac2{G(1)}
|\nabla v^0|(\Omega)$. This ends the proof of Step 6. Consequently
we get
$$II^0\geq J^0(u^0)\quad \mbox{with}\quad u^0\in B.$$

\noindent \textbf{Step 7: $\inf_{w\in B}J^0(w)>0$.}\\
First notice that $I^0=\inf_{w\in B}J^0(w) <+\infty$, because
$J^0(u^0)\leq
II^0\leq C_0<+\infty$.\\
Let us assume that $I^0=0$. Then we can consider a minimizing
sequence $w^k\in B$ such that $J^0(w^k)\longrightarrow 0$. By
definition of the $BV$-norm, of $B$ and of $J^0$, we see that the
sequence $(w^k)_k$ is bounded in $BV(\Omega)$:
$$\|w^{k}\|_{BV(\Omega)}\le \|w^{k}\|_{L^1(\Omega)}+\|\nabla
w^{k}\|_{L^1(\Omega)}= \displaystyle{|\Omega|
+\frac{1}{c}J^0(w^k)\leq |\Omega| +\frac{C_0}{c}}.$$ From the
compactness result for $BV$ (proposition \ref{pro:11}), up to
extract a subsequence, we get the existence of a function
$w^\infty\in BV(\Omega)$, such that
\begin{equation}\label{eq:27}
w^{k} \longrightarrow w^\infty \quad \mbox{in}\quad L^1(\Omega)
\end{equation}
and
$$0=\liminf_{k\to +\infty} |\nabla w^{k}|(\Omega) \ge |\nabla
w^{\infty}|(\Omega)$$ Therefore $w^\infty$ is constant on
$\Omega$, and because of (\ref{eq:27}) and $\int_\Omega w^k=0$, we
have
$$\int_\Omega w^\infty =0$$
and then $w^\infty \equiv 0$ on $\Omega$.\\
On the other hand, because of (\ref{eq:27}), up to extract a
subsequence, we have
$$w^{k}(x) \longrightarrow w^\infty(x) \quad \mbox{a.e. in}\quad \Omega$$
and then $w^\infty(x)=\pm 1$ a.e. in $\Omega$. Contradiction.\\
Then $\inf_{w\in B}J^0(w)>0$. This ends the proof of proposition
\ref{pro:10}.$\hfill\square$


\bigskip
\bigskip
\bigskip

\noindent \textbf{ Acknowledgment}\\
The authors would like to thank R.V. Kohn for enlightening
discussions. They also thank the referee of the paper whose
remarks and comments helped them to improve the presentation of
the article and to correct a confusion mistake in the statement of
Theorem \ref{th:4}.


\noindent A. El Soufi\\
Laboratoire de Math\'ematiques et Physique Th\'eorique\\ UMR 6083
du CNRS\\ Universit\'e de Tours\\ Parc de Grandmont\\ F- 37200
Tours, France\\
elsoufi@univ-tours.fr

\vskip 1cm

\noindent M. Jazar\\
Mathematics department\\
Lebanese University\\
P.O. Box 155-012, Beirut Lebanon\\
mjazar@ul.edu.lb

\vskip 1cm

\noindent R. Monneau\\
CERMICS - Ecole Normal des Ponts et Chauss\'ees\\
6 et 8 avenue Blaise Pascal Cit\'e Descartes - Champs sur Marne\\
77455 Marne la Vall\'ee Cedex 2, FRANCE\\
monneau@cermics.enpc.fr

\end{document}